\documentclass[11pt]{article}
\usepackage[utf8]{inputenc}
\usepackage[english]{babel}
\usepackage[margin=2.5cm,top=3cm]{geometry}
\usepackage{amsmath,amssymb,dsfont}
\usepackage{color}
\usepackage{theorem}
\usepackage{caption}
\usepackage{tikz}
\usepackage[pagewise]{lineno}

\usepackage{graphicx} 
\usepackage{epstopdf}
\usepackage{float}

\usepackage[colorlinks=true]{hyperref}

\newtheorem{defi}{Definition}[section]
\newtheorem{teo}[defi]{Theorem}
\newtheorem{prop}[defi]{Proposition}
\newtheorem{lem}[defi]{Lemma}
\newtheorem{cor}[defi]{Corollary}

{\theorembodyfont \rmfamily \newtheorem{remark}[defi]{Remark}}

\newcommand{\dis}{\displaystyle}

\newcommand{\R}{\mathds{R}}

\newcommand{\jnt}{\dis\int}

\newcommand{\avec}{\mathbf{a}}
\newcommand{\xvec}{\mathbf{x}}
\newcommand{\uvec}{\mathbf{u}}

\newcommand{\va}{\varphi}

\newcommand{\black}{\color{black}}

\newcommand{\fin}{\hfill$\square$}

\begin{document}

\title{ \huge 
\textbf{ Inverse problems for one-dimensional fluid-solid interaction models}}

\author{
J. Apraiz\thanks{Universidad del Pa\'is Vasco, Facultad de Ciencia y 
Tecnolog\'ia, Dpto.\ Matem\'aticas, Barrio Sarriena s/n 48940 Leioa (Bizkaia), 
Spain. E-mail: {\tt jone.apraiz@ehu.eus}.}
\ \ A. Doubova\thanks{Universidad de Sevilla, Dpto.\ EDAN e IMUS, 
Campus Reina Mercedes, 41012~Sevilla, Spain. E-mail: {\tt doubova@us.es}.},
\ \ E. Fern\'andez-Cara\thanks{Universidad de Sevilla, Dpto.\ EDAN e IMUS, 
Campus Reina Mercedes, 41012~Sevilla, Spain. E-mail: {\tt cara@us.es}.}, 
\ \ M. Yamamoto\thanks{The University of Tokyo, Graduate School of 
Mathematical Sciences, Komaba Meguro 153-8914 Japan 
E-mail:  {\tt myama@next.odn.ne.jp}.}
}

\date{}

\maketitle

\begin{abstract}
   We consider a one-dimensional fluid-solid interaction model governed by the Burgers equation with a time varying interface.
   We discuss on the inverse problem of determining the shape of the interface from Dirichlet and Neumann data at one end point of the spatial interval.
   In particular, we establish uniqueness results and some conditional stability estimates.
   For the proofs, we use and adapt some lateral estimates that, in turn, rely on appropriate Carleman and interpolation inequalities.
\end{abstract}

\vspace*{0,5in}
 





\section{Introduction}\label{sec_introduction}

   We will consider a nonlinear system that models the interaction of a one-dimensional fluid evolving in~$(-1,1)$ and a solid particle. 
   It will be assumed that  the velocity of the fluid is governed by the viscous Burgers equation at both sides of the point mass location~$y=p(t)$. 
   For simplicity, {\black it will be accepted that the fluid density is constant and equal to~$1$} and the solid particle has unit mass.

   {\black
   For any~$p$ at least in~$C^0([0,T])$ satisfying~$|p(t)| \leq 1$ for all~$t \in [0,T]$, let us introduce the open sets
   \[
   \begin{array}{c}
Q(p) =\{(x,t)\in\R^2 : -1<x<1,\ \ x\neq p(t),\ \  0<t<T\}, \\
Q_\ell (p) = \{(x,t)\in Q(p) : p(t) > x\} \ \text{ and } \
Q_r (p) = \{(x,t)\in Q(p) : p(t) < x\}.
   \end{array}
   \]}
   
   {\black On the other hand, the \emph{jump} of the function~$f$ at the point~$x$ will be denoted in the sequel by~$[f](x)$, that is,
   \[
[f](x) := \lim_{s\to0^+}f(x+s)-\lim_{s\to 0^-}f(x+s).
   \]}

   We will consider fluid-particle systems of the form
\begin{equation} \label{e_1}
     \left\{
    \begin{array}{ll}
    w_t-w_{xx}+ww_x=0, & (x,t) \in Q(p), 
    \\[1mm]
    w(p(t),t)=p'(t), \quad [w_x](p(t),t) = p''(t), & t \in (0,T),     
    \\[1mm]
     w(-1,t)=\alpha(t), \quad w(1,t)=\eta(t), & t \in (0,T),     \\[1mm]    
    w(x,0)=w_0(x), & x\in (-1,1), \\[1mm]
    p(0)=q_0,\quad p'(0)=q_1, &
    \end{array}
    \right.
\end{equation}
where (at least)~$w_0\in L^2(-1,1)$, $\alpha,\eta\in C^0([0,T])$, {\black $|q_0| < 1$} and~$q_1\in\R$.

   Here, $w(x,t)$ is the velocity of the fluid particle located at~$x$ at time~$t$, $p(t)$ is the position occupied by the particle at time~$t$ and~$\alpha$ and~$\eta$ are Dirichlet data.
   {\black It is assumed that~$w_0$, $q_0$ and~$q_1$ are initial data respectively for the fluid velocity, the particle position and the particle velocity.}

   {\black
   The first condition at~$x=p(t)$ in~$(\ref{e_1})$ means that the velocity  of the fluid and the solid mass coincide at this point.
   
   In the second condition, we state {\it Newton's law:}
   the force exerted by the fluid on the particle equals the product of the particle mass and its acceleration.
   Thus, if we introduce the notation~$u:=w|_{Q_\ell (p)}$ and~$v:=w|_{Q_r (p)}$, the jump condition at the points
$(p(t),t)$ can be written in the form
   \begin{equation} \label{eq_jump}
(v_x - u_x)(p(t),t)=p''(t), \quad t \in (0,T).
   \end{equation}
   }

   {\black The previous system can be viewed as a preliminary simplified version of other 
more complicate and more realistic models in higher dimensions that we plan to analyze in the future.
   For example, it is meaningful to consider a system governed by the Navier-Stokes equations around a moving sphere that interacts with the fluid.
   
   More precisely, let~$R_0$, $R_1$ and~$R_{\rm ext}$ be given with $0 < R_0 < R_{\rm ext}$, let~$B(R_0)$ be the open ball centered at~$0$ of radius~$R_0$, let us set
   \[
\Sigma_{\rm ext} = \{ (\xvec,t) : |\xvec| = R_{\rm ext}, \ \ 0 < t < T \}
   \]
and, for any $R \in C^0([0,T])$ with $R_0 \leq R(t) \leq R_{\rm ext}$, let us introduce
   \[
   \begin{array}{c}
Q(R) =\{ (\xvec,t) : R(t) < |\xvec| < R_{\rm ext}, \ \  0 < t < T \}, \\
\Sigma(R) =\{ (\xvec,t) : |\xvec| = R(t), \ \  0 < t < T \}.
   \end{array}
   \]
   Then, it makes sense to search for functions $\uvec$, $p$ and $R$ satisfying
   \[
    \left\{
    \begin{array}{ll} \dis
    \uvec_t - \nu \Delta \uvec + (\uvec\cdot\nabla)\uvec + \nabla p = 0,
    \ \ \nabla \cdot \uvec = 0, & (\xvec,t) \in Q(R), 
    \\[1mm]  \dis
    \uvec(\xvec,t) = \avec(\xvec,t), & (\xvec,t) \in \Sigma_{\rm ext} ,     
    \\[1mm]  \dis
    \uvec(\xvec,t) = \frac{\dot{R}(t)}{R(t)}\xvec,
    & (\xvec,t) \in \Sigma(R),     
    \\[1mm]  \dis
    \uvec(\xvec,0) = \uvec_0(\xvec), & \xvec\in B(R_0),
    \\[1mm]  \dis
    R(0) = R_0, \ \ \dot{R}(0) = R_1, &
    \\[1mm]  \dis
    \left(-p \mathbf{Id.} + 2\nu D\uvec \right) \frac{\xvec}{R(t)} 
    = \ddot{R}(t)\,\xvec,
    & (\xvec,t) \in \Sigma(R) ,
    \end{array}\right.
   \]
where the constant~$\nu > 0$ and the radially symmetric fields~$\avec = \avec(\xvec,t)$ and~$\uvec_0 = \uvec_0(\xvec)$ are given, $\bf{Id.}$ is the identity matrix and~$D\uvec$ denotes the symmetrized gradient of~$\uvec$, that is,
   \[
D\uvec = \frac{1}{2}(\nabla\uvec + \nabla\uvec^T).
   \]

   A related question is whether we can determine the function $R$ from~$\avec$, $\uvec_0$ and exterior boundary observations
   \[
F = \left(-p \mathbf{Id.} + \nu D\uvec \right) \frac{\xvec}{R_{\rm ext}}
   \]
on (a part of) $\Sigma_{\rm ext}$.
   }
   
   This justifies the relevance of the analysis of inverse problems for~\eqref{e_1}.

\
   
   As far as we know, the first works where the simplified model~\eqref{e_1} has been considered are~\cite{VZ2} and~\cite{VZ1}.
   There, the authors allowed the spatial variable to take any value in~$\R$ instead of~$(-1,1)$.
   In particular, in~\cite{VZ2}, the authors proved the existence and uniqueness of a solution and described its large-time behavior for just one solid mass submerged in the fluid.
   In~\cite{VZ1}, similar result were established in the case of various rigid bodies immersed in the fluid.
   
   {\black
   These results were later extended to a multi-dimensional framework in~\cite{MZ}. 
   Let us also mention that  the controllability properties of a system 
similar to~\eqref{e_1} have been analyzed in~\cite{{DF2}} and~\cite{LTT}.}

\

   In what concerns the {\it direct problem,} that is, to find appropriate~$w$ and~$p$ verifying the equation and additional conditions in~\eqref{e_1}, it can be shown that, if~$\vert q_0\vert + \vert q_1 \vert + \Vert w_0\Vert_{H^1_0(-1,1)}$ 
is sufficiently small, there exists a solution~$(w,p)$ to~\eqref{e_1} 
with~$w \in  C^0([0,T];H^1(-1,1))$, $w_{xx} \in L^2(0,T;L^2(-1,1))$ {\black and~$p\in H^2(0,T)$;}
%
   see for example~\cite[Theorem 1.1]{LTT}.
   
   {\black In fact, the result in~\cite{LTT} only states that~$p \in C^1([0,T])$.
   However, the regularity of the restrictions of~$w$ to~$Q_\ell(p)$ and~$Q_r(p)$ 
shows that the a.e.\ defined function $t \mapsto [w_x](p(t),t)$ is square-integrable and, consequently, $p'' \in L^2(0,T)$.}

   The {\it inverse problem} related to system~\eqref{e_1} we are interested in is 
the following:

\

\noindent
\textbf{Inverse problem -} Given the data~$T>0$, $q_0\in (-1,1)$, $q_1\in\R$ and~$\alpha \in C^0([0,T])$ and the observation~$\beta$ with $\beta(t) = w_x(-1,t)$ for~$t\in (0,T)$, find~$\eta:= w(1,\cdot)$.

\

   In this paper, we will study related uniqueness and stability properties. 
   In particular, we will give answers to questions like the following: 

\

{\black
\noindent
\textbf{Global uniqueness -} Let~$(w_i, p_i)$ be a solution to~\eqref{e_1} associated to some $T$, $q_0$, $q_1$ and~$\alpha$ for~$i = 1,2$.
   Assume that  the corresponding observations coincide at~$x=-1$, that is, $w_{1,x}(-1,t) =  w_{2,x}(-1,t)$ for~$0<T_1<t<T_2<T$.
   Then, do we have~$p_1=p_2$  and~$w_1=w_2$?

\

\noindent
\textbf{Global stability -} Let~$(w_i, p_i)$ be as before and set~$\beta_i := w_{i,x}(-1,\cdot)$ and~$\eta_i = w_i(1,\cdot)$ for~$i = 1,2$.
   Is there any estimate of the kind
   \[
\|\eta_1 - \eta_2\|_{L^\infty(T_1,T_2)} + \|p_1 - p_2\|_{L^\infty(T_1,T_2)} \leq 
\phi(\|\beta_1 - \beta_2\|_{L^\infty(0,T)})
   \]
for some continuous function $\phi : \R_+ \mapsto \R_+$ satisfying $\lim_{s\to0^+} \phi(s) = 0$?
}

\

   {\black The paper is organized as follows.
   
   First, in Section~\ref{sec_preliminary}, we prove a preliminary fundamental lemma that plays a key role in the proof of conditional stability.
   It provides estimates of the traces on the interface~$x = p(t)$ of the difference of two solutions to~\eqref{e_1} in terms of the boundary data and observations.}
   
   In Section~\ref{sec_unilateral}, we establish a stability estimate and then the uniqueness of the lateral inverse problem corresponding to the system satisfied in the left part~$Q_\ell(p)$ of the whole domain.
   By reflection, similar results are fulfilled by the solution to the system satisfied in~$Q_r(p)$.
   
   Section~\ref{sec_global} is devoted to establish a global stability and uniqueness result for the inverse problem in the whole domain~$Q(p)$.
   

\

\section{Preliminaries}\label{sec_preliminary}

   As already said, the main result in this section is crucial for the proof of a local stability property that will be established in~Section~\ref{sec_unilateral} (see~Proposition~\ref{prop_unilateral}).


\

\begin{lem}\label{lemma_key} 
   Let us assume that
   \begin{equation} \label{eq_2}
   \begin{cases}
\ u_t -u_{xx} + au_x + bu = 0, &   (x,t)\in Q_\ell (p), \\
\ u(-1,t) = \alpha(t), \quad  u_x(-1,t) = \beta(t), & t\in (0,T),
   \end{cases}
   \end{equation}
with~$a, b\in L^{\infty}(Q_{\ell}(p))$, $u\in H^2(Q_{\ell}(p))$ and
there exist constants~$M>0$ and~$\delta \in (0,1)$ such that  
   \begin{equation}\label{e_2}
\Vert u\Vert_{H^2(Q_{\ell}(p))}\leq M, \ \  {\black \Vert p\Vert_{H^2(0,T)}\leq M \ \text{ and }} \ |p(t)|\leq 1-\delta \quad \forall\, t\in[0,T].
   \end{equation}
   Then:
   
\begin{itemize}
   
\item[a)] For any~$\epsilon>0$, there exist constants~$K_{\epsilon}>0$ 
and~$\theta_{\epsilon}\in (0,1)$ such that 
   \begin{equation}\label{e_4a}
|u(p(t),t)|\leq\frac{K_{\epsilon}}{\left|\log{\frac{1}{k}}\right|
^{\theta_{\epsilon}}} \quad \forall\, t\in[\epsilon,T],
   \end{equation}
provided~$\alpha$, $\beta$ and~$k$ satisfy 
   \begin{equation}\label{e_3}
0\le \Vert \alpha\Vert_{L^2(0,T)}+\Vert \beta\Vert_{L^2(0,T)}<k<1.
   \end{equation}

\item[b)] In particular, if~$\alpha\equiv 0$ and~$\beta\equiv 0$ in~$(0,T)$, 
then~$u\equiv 0$ in~$Q_{\ell}(p)$.

\end{itemize}
\end{lem}


%

\

\noindent 
\textbf{Proof:}  

   The proof of part a) can be obtained by adapting some arguments in~\cite{Ya} that rely on appropriate Carleman estimates.
   Carleman estimates were first used in~\cite{BK-1} to establish uniqueness and stability results for inverse problems;
   see also~\cite{K-1, K-2}.   

   {\black The argument is decomposed into three steps.}

\

\noindent
$\bullet$ \textbf{Step 1:} First, we introduce the change of variables
   \[
y = \frac{x+1}{p(t)+1}, \quad z(y,t) = u(x,t)\quad\text{for}\quad 
(x,t)\in Q_{\ell}(p).
   \]

   In this way, the values of~$y$ remain in the new spatial domain~$(0,1)$, $Q_{\ell}(p)$ is transformed into~$Q_\ell :=(0,1)\times(0,T)$ and~$z$ satisfies the system
   \[
   \begin{cases}
\ z_t -\dfrac{1}{(p(t)+1)^2}\,z_{yy} + \dfrac{p'(t)}{p(t)+1} A(y,t) z_y 
+ \dfrac{1}{p(t)+1} B(y,t) z = 0, &   (y,t)\in Q_\ell, 
\\[4mm]
\ z(0,t) = \alpha(t), \quad  z_y(0,t) = (p(t)+1) \beta(t), & t\in (0,T),
   \end{cases}
   \]
where
   \[
A(y,t) :=  a\left(\frac{x+1}{p(t)+1}\right) \quad\text{and}\quad B(y,t) :=  b\left(\frac{x+1}{p(t)+1}\right).
   \] 

   Then, we perform a second change of variables:
   \[
x^*=2-2y, \quad u^*(x^*,t) = z(y,t)\quad\text{for}\quad (y,t)\in Q_{\ell}
   \]
and we now have
   \begin{equation}\label{eq.ustar}
\begin{cases}
\ u^*_t -\dfrac{1}{(p(t)+1)^2}\,u^*_{x^*x^*} 
- \dfrac{2p'(t)}{p(t)+1} A^*(x^*,t)u^*_{x^*} 
+ \dfrac{1}{p(t)+1} B^*(x^*,t)u^* = 0, &   (x^*,t)\in Q^*, 
\\[5mm]
\ u^*(2,t) = \alpha(t), \quad  u^*_{x^*}(2,t) = -\dfrac{1}{2}(p(t)+1) \beta(t), 
& t\in (0,T),
\end{cases}
   \end{equation}
where~$Q^*:=(0,2)\times(0,T)$.

   After skipping the stars in these variables, coefficients and sets, we see that  the task is reduced to prove the existence of~$K_{\epsilon}$ and  $\theta_{\epsilon}$ such that  
   \begin{equation}\label{e_6}
|u(0,t)|\leq\frac{K_{\epsilon}}{\left|\log{\frac{1}{k}}\right|
^{\theta_{\epsilon}}} \quad \forall\, t\in[\epsilon,T].
   \end{equation}
   
   Thus, the rest of the proof is devoted to establish~\eqref{e_6}.


\

\medskip\noindent
$\bullet$ \textbf{Step 2:} Let us start with the proof of the following intermediate estimate: 
   \begin{equation}\label{e_7}
|u(0,\overline{t})|\leq\frac{K_{0,\epsilon}}{\left|\log{\frac{1}{F_{\epsilon}}}
\right|^{\theta_0}} \quad \forall\, \overline{t}\in[\epsilon,T],
   \end{equation}
where~$\theta_0\in (0,1)$ is independent of~$\epsilon$ and  
   \[
F_{\epsilon}:=\sup_{x\in[1,2]}\left(\Vert u(x,\cdot)\Vert 
_{L^2(\epsilon, T)}+\Vert u_x(x,\cdot)\Vert_{L^2(\epsilon, T)}\right) .
   \]
 
   To this purpose, let us fix~$\bar{t}\in[\epsilon,T]$ and let us introduce a new change of variables: 
   \begin{equation}\label{eq.changev}
 \hat{t} = \dfrac{T}{2 \bar{t} }\, t,\quad  \hat{u}(x,\hat{t}) = u(x,t) \quad\text{for}\quad (x,t) \in (0,2)\times(0,T).
   \end{equation}
   Then, $\hat{u}$ is defined in~$\hat{Q} := (0,2)\times (0,T^2/(2\overline{t}))$ and satisfies a system similar to~\eqref{eq.ustar} with coefficients~$\hat{D}$, $\hat{A}$ and~$\hat{B}$ that are uniformly bounded for~$\overline{t} \ge \epsilon >0$:
   \begin{equation}\label{eq.uwidehat}
   \begin{cases}
\ \hat{u}_{\hat{t}} -\hat{D}(\hat{t})\hat{u}_{xx} 
- \hat{A}(x,\hat{t})\hat{u}_{x} + \hat{B}(x,\hat{t})u = 0, 
& (x,\hat{t})
\in \hat{Q}, 
\\[2mm]
\ \hat{u}(2,\hat{t}) = \hat{\alpha}(\hat{t}), \quad  
\hat{u}_x(2,\hat{t}) 
= \hat{\beta}(\hat{t}), & \hat{t}\in (0,T^2/(2\overline{t})).
   \end{cases}
   \end{equation}
   Therefore, what we have to prove is~\eqref{e_7} for~$\hat{t}=T/2$, that is, 
   \begin{equation}\label{uint2}
|\hat{u}\left(0,T/2\right)| \leq \frac{K_{0,\epsilon}}
{\left|\log{\frac{1}{F_{\epsilon}}}\right|^{\theta_0}}.
   \end{equation}

   In the sequel, we will distinguish two cases.

\

\begin{figure}[h!]
\centering
\includegraphics[width = 7cm]{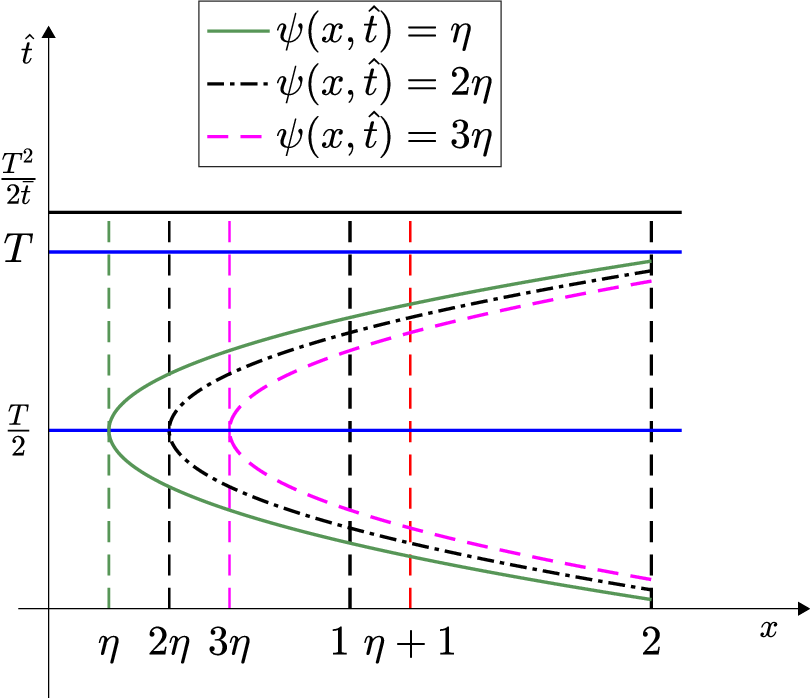}
\caption{The curves~$\psi = \eta$, $\psi = 2\eta$ and~$\psi = 3\eta$. \label{fig.Qetas}} 
\end{figure}

\begin{figure}[h!]
\centering
\includegraphics[width = 13cm]{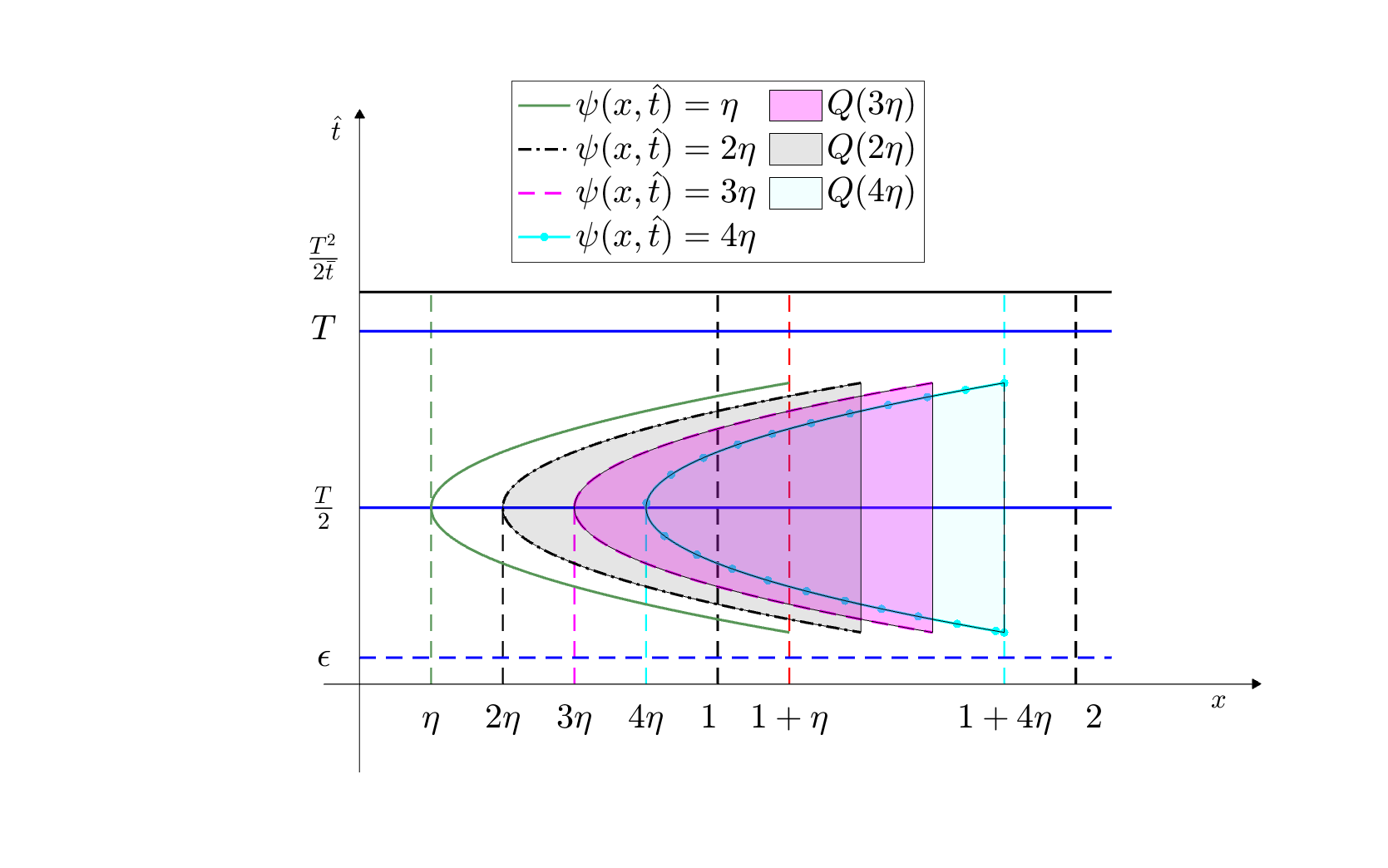}
\vspace{-0.8cm}
\caption{The sets~$Q(2\eta)$, $Q(3\eta)$ and~$Q(4\eta)$.} 
\label{fig.QetasBIS}
\end{figure}

\noindent
{\black\textit{Case 1:} $\overline{t} \le T/2$.
   
   Let us introduce 
   \begin{equation}\label{eq.psi}
\psi(x,\hat{t}) := x-\rho | \hat{t} - T/2 |^2
\quad\text{and}\quad \varphi(x,\hat{t}) := e^{\lambda \psi(x,\hat{t})},
   \end{equation}
where~$4/(T-2\epsilon)^2 < \rho < 16/T^2$ and $\lambda>0$ is sufficiently large and the sets
   \[
Q(\eta) := \{(x,\hat{t})\in(0,2)\times(0,T): 0<x<1+\eta,\, \psi(x,\hat{t})>\eta\}
   \]
for $\eta \in (0,1)$;
   see Figure~\ref{fig.Qeta}.
   
\begin{figure}[h!]
\centering
\includegraphics[width = 11cm]{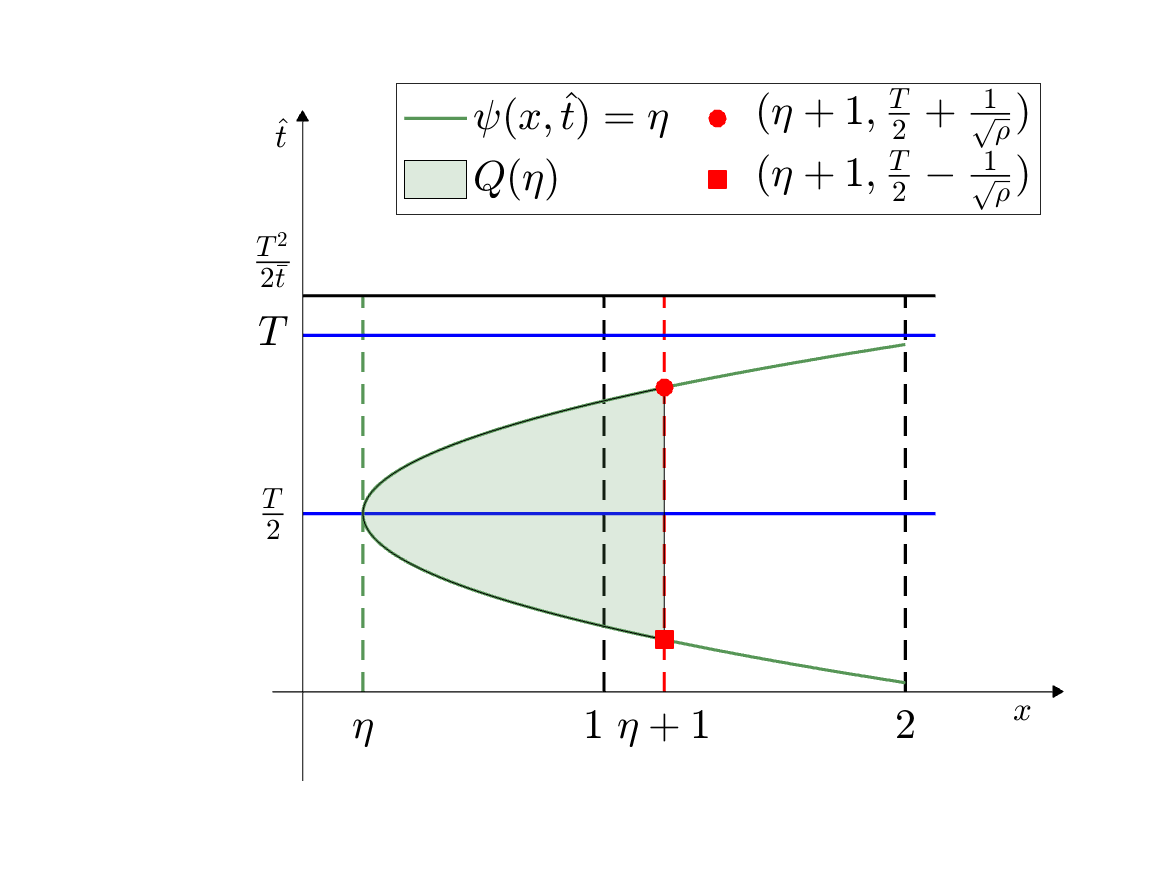}
\vspace{-0.8cm}
\caption{The set~$Q(\eta)$ when~$\bar{t} \leq T/2$. \label{fig.Qeta}} 
\end{figure}

   For any~$\eta \in (0,1/4)$, let us introduce a function~$\chi_{\eta} \in C^{\infty}(\R^2)$ satisfying~$\chi_{\eta}(x,\hat{t}) = 1$ for~$\psi(x,\hat{t})\geq 3\eta$, $\chi_{\eta}(x,\hat{t}) = 0$ for~$\psi(x,\hat{t})\leq 2\eta$ and, for instance, $0\le \chi_{\eta}\le 1$. 
   Then
\begin{equation}\label{eq.chi}
 \chi_{\eta}(x,\hat{t}) = 
 \begin{cases}
\ 1 & \text{ in } Q(3\eta) ,\\[1mm]
\  0 & \text{ in } Q{(\eta)}\setminus Q(2\eta) ;
\end{cases}
\end{equation}
   see Figures~\ref{fig.Qetas} and~\ref{fig.QetasBIS}.}

   Moreover, 
   \begin{equation}\label{eq.estxi}
| \chi_{\eta,x} | + | \chi_{\eta,\hat{t}} |  \leq  \dfrac{C}{\eta}
\quad\text{and}\quad 
| \chi_{\eta,xx} | + | \chi_{\eta,x\hat{t}} | + | \chi_{\eta,\hat{t}\hat{t}} |  
\leq \dfrac{C}{\eta^2}.
   \end{equation}

   Then, let us set~$v_{\eta}:=\hat{u}\chi_{\eta}$.
   At this point, we can use a global Carleman estimate for~$v_{\eta}$ in~$Q(\eta)$ (see~\cite[Theorem 3.2]{Ya}).
   
   {\black Thus, let us rewrite the PDE in~\eqref{eq.uwidehat} in the form~$E \hat{u} = 0$.}
   There exist~$s_0,\lambda_0 > 0$ such that, for any~$s \geq s_0$ and any~$\lambda \geq \lambda_0$, one has
   \begin{equation}\label{e_8}
s^3\lambda^4\iint_{Q(\eta)}e^{2s\va}\va^3|v_{\eta}|^2 \, dx\,dt
\leq C_{0,\epsilon}\iint_{Q(\eta)}e^{2s\va}|H|^2 \, dx\,dt
+C_{1,\epsilon} e^{Cs}G^2_{\sigma} \,,
   \end{equation}
where~$C_{0,\epsilon}$ and~$C_{1,\epsilon}$ are constants steming from the uniform bounds of the coefficients in~\eqref{eq.uwidehat} for all~$\hat{t} \geq \epsilon > 0$ and~$H := E v_\eta$.
   Note that $H$ can be nonzero only in~$Q(2\eta)\setminus Q(3\eta)$.

   In~\eqref{e_8}, we can take 
   \[
\displaystyle{G^2_{\sigma}
:=\sup_{x\in[1,2]}\left(\Vert \hat{u}(x,\cdot)\Vert^2_{L^2(\sigma, T-\sigma)}
+\Vert \hat{u}_x(x,\cdot)\Vert^2_{L^2(\sigma, T-\sigma)}\right)} ,
   \]
with~$\sigma = T/2 - 1/\sqrt{\rho}$.
   Then, thanks to the choice of~$\rho$ in~\eqref{eq.psi} one has~$[\sigma,T-\sigma]\subset (\epsilon, T)$.

   Without loss of generality, we can assume that~$G_{\sigma}<1$. 
 
   {\black In the sequel, we will denote by~$C$ (resp.~$C_\epsilon$) a generic positive constant independent of~$s$ and~$\epsilon$ (resp.~independent of~$s$).
   Of course, both constants can change from line to line.

   Taking into account the definition of~$Q(\eta)$  (see Figure~\ref{fig.QetasBIS}), we have from~\eqref{e_8} that
   \begin{equation}\label{e_8bis}
s^3\lambda^4\iint_{Q(4\eta)\cap Q(\eta) }e^{2s\va}\va^3|\hat{u}|^2\, dx\,dt
\leq C_{\epsilon}\iint_{Q(2\eta)\setminus Q(3\eta)}e^{2s\va}|H|^2\, dx\,dt
+ C_{\epsilon} e^{Cs}G^2_{\sigma}
   \end{equation}
for any~$s \geq s_0$ and any~$\lambda\geq\lambda_0$. 

   Now, in view of~\eqref{eq.estxi}, using the fact that~$\varphi\ge \mu_2 := e^{4\lambda \eta}$ in~$Q(4\eta)\cap Q(\eta)$ and~$\varphi\le \mu_1:= e^{3\lambda \eta}$ in~$Q(2\eta)\setminus Q(3\eta)$, 
we obtain from~\eqref{e_8bis} that  
   \begin{align*}
e^{2s\mu_2}\iint_{Q(4\eta)\cap Q(\eta)}|\hat{u}|^2\, dx\,dt
&\leq 
\dfrac{C_\epsilon}{\eta^4} e^{2s\mu_1}
\iint_{Q(2\eta)\setminus Q(3\eta)} \left(|\hat{u}|^2+ |\hat{u}_x|^2\right) \, dx\,dt 
+ C_{\epsilon} { e^{Cs}}G_\sigma^2
   \end{align*}
and, therefore,  
   \begin{align*}
\iint_{Q(4\eta)\cap Q(\eta)}|\hat{u}|^2\, dx\,dt
&\leq 
\dfrac{C_\epsilon}{\eta^4} e^{-2s(\mu_2-\mu_1)}\Vert  \hat{u}\Vert_{H^1(Q)} ^2
+ C_{\epsilon} { e^{Cs}} G_\sigma^2 \quad \forall s\geq s_0.
   \end{align*}

   From the definition of~$G_\sigma$ and assmption~\eqref{e_2}, one has:
   \begin{align*}
\Vert \hat{u}\Vert^2_{L^2(Q(4\eta))}&\leq 2\Vert \hat{u} \Vert^2
_{L^2(Q(4\eta)\cap Q(\eta))}+2\Vert \hat{u}\Vert^2
_{L^2(Q(4\eta)\setminus Q(\eta))}
\\[2mm]
&\leq\frac{C_\epsilon}{\eta^4}e^{-2s(\mu_2-\mu_1)}\Vert \hat{u}\Vert^2_{H^{1}(Q)}
+C_{\epsilon} { e^{Cs}}G_\sigma^2
\\[2mm]
&\leq\frac{C_\epsilon}{\eta^4}e^{-2s\lambda \eta}M^2+C_{\epsilon} {e^{Cs}}G_\sigma^2 .
   \end{align*}
   Thus, if we fix~$\lambda$, we take~$s'=s+s_0$ and then rename~$s'$ as~$s$,  we deduce that  
   \begin{equation}\label{ineq-new}
\Vert \hat{u}\Vert_{L^2(Q(4\eta))}
\leq
C_\epsilon\left(\dfrac{1}{\eta^2}e^{-s\lambda\eta}M + G_\sigma{e^{Cs}}\right)
\quad\forall s\geq 0.
   \end{equation}

   
   Now, in order to get the best estimate, we minimize this quantity by choosing
   \[
s = \frac{1}{C + \lambda\eta} \, \log\left(\frac{M\lambda}{CG_\sigma\eta}\right)
   \]
and we deduce that
   \begin{equation*}
\Vert \hat{u}\Vert_{L^2(Q(4\eta))} \leq C_\epsilon \,
\frac{M^{C/(C + \lambda\eta)}}{\eta^{1 + C/(C + \lambda\eta)}} \, G_{\sigma}^{\lambda\eta/(C+\lambda\eta)}
\quad\text{for}\quad 0<\eta<\frac{1}{4} .
   \end{equation*}

   Then, using the classical Sobolev embedding and an interpolation inequality, for any~$\ell \in (1,2)$, we obtain
   \begin{align*}
\Vert \hat{u}\Vert_{L^{\infty}(Q(4\eta))}
&\leq 
C\Vert \hat{u}\Vert_{H^\ell(Q(4\eta))}
\leq 
C\Vert \hat{u} \Vert^{\ell/(\ell+1)}_{H^{\ell+1}(Q(4\eta))} \Vert \hat{u} \Vert^{1/(\ell+1)}_{L^2(Q(4\eta))}
\\[2mm]
&\leq
C_\epsilon M^{\ell/(\ell+1)} \left(\frac{M^{C/(C + \lambda\eta)}}{\eta^{1 + C/(C + \lambda\eta)}}\right)^{1/(\ell+1)}
G_{\sigma}^{\lambda\eta/((C+\lambda\eta)(\ell+1))}
   \end{align*}
whenever $0 < \eta < 1/4$.
   In particular, after a choice of~$\ell$ in~$(1,2)$, one has for some $a \in (1,2)$ independent of~$\eta$ that
   \[
| \hat{u}(4\eta,T/2) | \leq \frac{C_\epsilon}{\eta^{1-a/2}}MG_{\sigma}^{C\eta} \quad \forall \eta \in (0,1/4).
   \]
   After integration with respect to~$\eta$ and a variable change~$x=4\eta$, we obtain:
   \begin{align*}
\int_0^1| \hat{u} (x,T/2)|^2\,dx
= & 
4\int_0^{1/4}\Big | \hat{u} (4\eta,T/2)|^2\, d\eta
\leq
C_\epsilon M^2 \jnt_0^{+\infty}\frac{1}{\eta^{2-a}}G_{\sigma}^{2C\eta}\, d\eta
\\[3mm]
= & C_\epsilon M^2 \jnt_0^{+\infty}\frac{1}{\eta^{2-a}}
e^{-2C\log{(1/G_\sigma)}\eta}\, d\eta
= C_\epsilon M^2 (\log{(1/G_\sigma)})^{-(a-1)}.
   \end{align*}
   }
   
   Thus, using again interpolation, we see that, for any $m \in (1/2,1)$,
   \begin{align*}
\Vert \hat{u}\left(\cdot,T/2\right) \Vert
_{L^{\infty}(0,1)}
&\leq
\Vert \hat{u} \left(\cdot,T/2\right) \Vert
_{H^m(0,1)}
\leq
\Vert \hat{u}\left(\cdot,T/2\right) \Vert^m
_{H^1(0,1)} \Vert \hat{u} \left(\cdot,T/2\right) \Vert^{1-m}_{L^2(0,1)}
\\[4mm]
&
\leq C M^m \Vert \hat{u} \left(\cdot,T/2\right) \Vert^{1-m}_{L^2(0,1)}
\leq C_\epsilon M\frac{1}{\left| \log{(1/G_\sigma)} \right|^{\theta_0}},
   \end{align*}
where~$\theta_0=\frac{1}{2}(a-1)(1-m)$.

   We note that~$\theta$ can be taken arbitrarily close to~$0$ and is independent of~$\epsilon$.
   It is clear that~$G_\sigma \le F_\epsilon$. 
   Therefore, if~$\bar{t} \leq T/2$, we have~\eqref{uint2}, which gives~\eqref{e_7}.

\

\noindent
{\black\textit{Case 2:} $\bar{t} > T/2$.}

\begin{figure}[h!]
 \centering
 \includegraphics[width = 11cm]{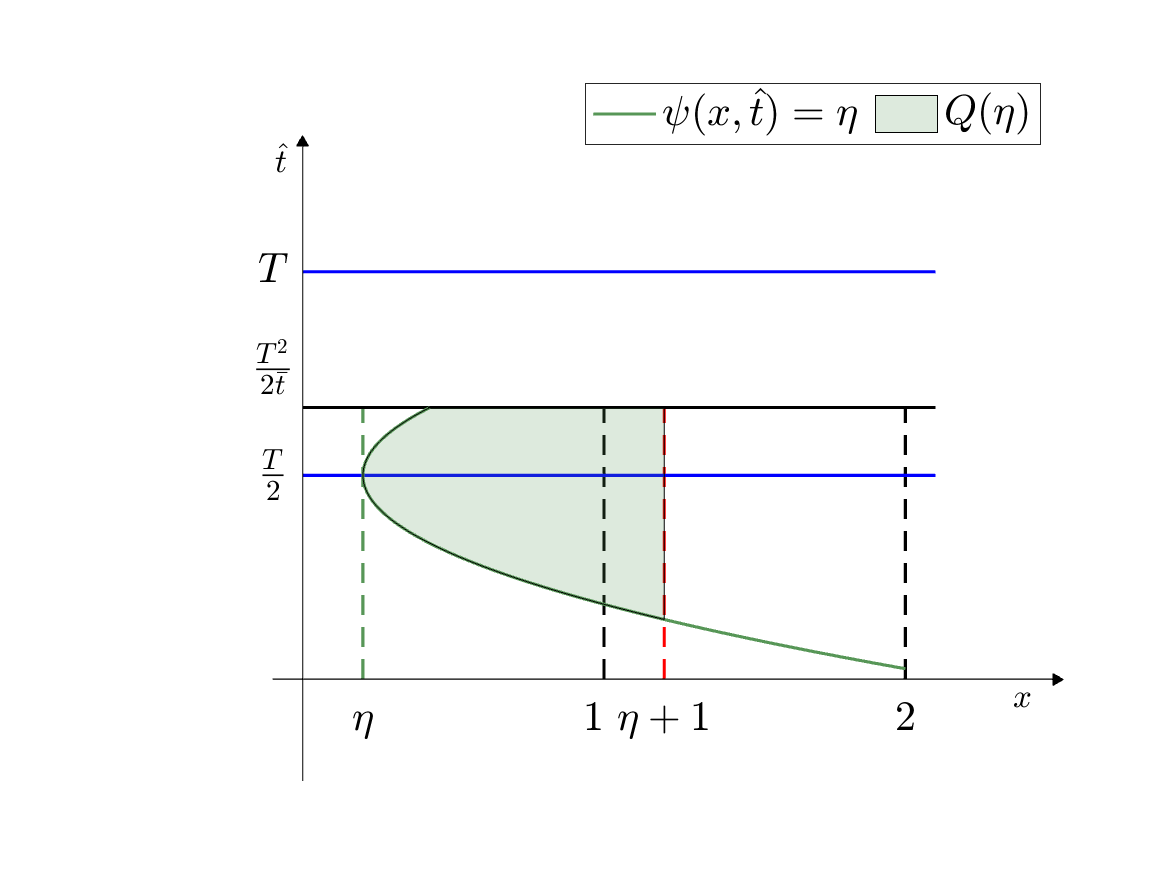}
\vspace{-0.8cm}
 \caption{The set~$Q(\eta)$ when~$\bar{t} > T/2$.\label{fig.QetaCase2}} 
 \end{figure}

   Let us introduce again the function~$\psi$, given by~\eqref{eq.psi}.
   In this case, we consider the sets
   \[
Q(\eta)=\{(x,\hat{t})\in(0,2)\times(0,T^2/(2\bar{t})):  0<x<1+\eta,\, \psi(x,t)>\eta\}
   \] 
for $0 < \eta < 1/4$; 
   see Figure~\ref{fig.QetaCase2}.
 
   Using once more~$v_{\eta} = \hat{u}\chi_{\eta}$ and  arguing as before, we get 
   \begin{equation}\label{e_8bis2}
s^3\lambda^4\iint_{Q(4\eta)\cap Q(\eta) }e^{2s\va}\va^3|\hat{u}|^2\, dx\,dt
\leq 
C_{\epsilon}\iint_{Q(2\eta)\setminus Q(3\eta)}e^{2s\va}|H|^2\, dx\,dt
+ C_{\epsilon} e^{Cs}G^2_{\sigma},
   \end{equation}
where we have set
   \[
\displaystyle{G^2_{\sigma}
:=\sup_{x\in[1,2]}\left(\Vert \hat{u}(x,\cdot)\Vert^2_{L^2(\sigma, T^2/(2\bar{t}))}
+\Vert \hat{u}_x(x,\cdot)\Vert^2_{L^2(\sigma, T^2/(2\bar{t}))}\right)}.
   \] 

   It is straightforward to see that the arguments applied in the previous case yield again~\eqref{uint2}.
   Consequently, we also have~\eqref{e_7} in this case.

\

\noindent
$\bullet$ \textbf{Step 3:} Let us establish the following lateral estimate of~$F_{\epsilon}$:
   for every~$\epsilon>0$, there exist constants~$C_{\epsilon}>0$ and~$\theta_{\epsilon}\in(0,1)$ such that 
   \begin{equation}\label{e_9}
G_\sigma\le F_{\epsilon}\leq C_{\epsilon}
M^{1-\theta_{\epsilon}}\left(\Vert \alpha\Vert_{L^2(0,T)}
+ \Vert \beta\Vert_{L^2(0,T)}\right)^{\theta_{\epsilon}}
+ C_{\epsilon}\left(\Vert \alpha\Vert_{L^2(0,T)}+\Vert \beta\Vert_{L^2(0,T)}\right).
   \end{equation}
 
   To this purpose, we first use Theorem 5.1 in~\cite{Ya} and deduce that
   \begin{equation}\label{e_10}
F_{\epsilon}\leq C_{0,\epsilon}M^{1-\theta_{0,\epsilon}}
\left(\Vert \alpha\Vert_{H^1(0,T)}+\Vert \beta\Vert_{L^2(0,T)}\right)^{\theta_{0,\epsilon}}
+ C_{\epsilon}\left(\Vert \alpha\Vert_{H^1(0,T)}+\Vert \beta\Vert_{L^2(0,T)}\right),
   \end{equation}
where~$\theta_{0,\epsilon} \in (0,1)$ and~$C_{\epsilon}>0$.

   Then, we use an interpolation inequality
   \begin{equation}\label{e_11}
\Vert \alpha\Vert_{H^1(0,T)}\leq\Vert \alpha\Vert^{1/2}_{H^2(0,T)}\Vert \alpha\Vert^{1/2}
_{L^2(0,T)}
\leq C\Vert u\Vert^{1/2}_{H^2(Q)}\Vert \alpha\Vert^{1/2}_{L^2(0,T)}
\leq C\ M^{1/2} \Vert \alpha\Vert^{1/2}_{L^2(0,T)}
   \end{equation}
and we observe that~\eqref{e_10} and~\eqref{e_11} imply~\eqref{e_9} when~$\Vert \alpha\Vert_{L^2(0,T)}+\Vert \beta\Vert_{L^2(0,T)}< k <1$.

   Finally, part b) of the lemma is straightforward:
   it suffices to use~\eqref{e_4a}  with~$\alpha=0$, $\beta = 0$ and~$k$ arbitrarily small.
\fin

\

\begin{remark}
   A similar result can be obtained when the right-hand side of the first equation in~\eqref{eq_2}, is not zero.
   Thus, if
   \[
u_t-u_{xx}+au_x+bu=f \ \text{ for } \ (x,t)\in Q_{\ell}(p)
   \]
and (for instance) $0 \leq \Vert f\Vert_{L^{\infty}(Q_{\ell}(p))}+\Vert \alpha\Vert_{L^2(0,T)}+\Vert \beta\Vert_{L^2(0,T)} < m < 1$, we can prove that, for any~$\epsilon>0$,
   \begin{equation}\label{e_4b}
|u(p(t),t)|\leq
\dfrac{K_{\epsilon}}{(\log{\frac{1}{m}})^{\theta_{\epsilon}}}
\quad\forall\, t\in[\epsilon, T]
   \end{equation}
for some~$K_{\epsilon}>0$ and some~$\theta_{\epsilon}\in (0,1)$.
\fin
\end{remark}

\

{\black\begin{remark}
   A more involved argument is needed to take~$\epsilon = 0$ in~\eqref{e_4a} and~\eqref{e_4b}.
   The stability rate is expected to be weaker than single logarithmic.
   This will be analyzed in a forthcoming paper.
\fin
\end{remark}}

\

\section{Lateral estimates and uniqueness}\label{sec_unilateral}
   
   This section is devoted to study the stability and uniqueness of~\eqref{e_1} on the left part of the domain, $Q_\ell (p)$.
   Later, we will extend these results to~$Q_r (p)$ and will obtain similar results in the whole domain~$Q(p)$. 

   Assume that 
   \begin{equation} \label{e_12}
\left\{
\begin{array}{ll}
u^i_t-u^i_{xx}+u^iu^i_x=0, & (x,t) \in Q_\ell (p_i),  \\[1mm]
u^i(-1,t)=\alpha^i(t), \quad    u^i_x(-1,t)=\beta^i(t),  & t \in (0,T), 
\\[1mm]
u^i(p_i(t),t)=p_{i}'(t), & t\in (0,T),      
\end{array}\right.
    \end{equation}
for~$i = 1, 2$.
   Let us formulate an inverse problem concerning the left part of the domain:

\

\noindent
{\black\textbf{Lateral uniqueness in~$Q_\ell (p)$:}
   Let~$(u^i, p_i)$, $i=1,2$ be two solutions to~\eqref{e_12} in~$Q_\ell (p_i)$.
   Assume that  the corresponding  observations coincide at the boundary~$x=-1$, that is,  
   \[
u_x^1(-1,t) =  u_x^2(-1,t) \ \text{ in some time interval~$(T_1,T_2)$.}
   \]
   Then, do we have~$p_1=p_2$ in~$(0,T)$ and~$u^1=u^2$ in~$Q_\ell (p)$ with~$p = p_1=p_2$?}

\

   {\black As before, we will denote in the sequel by~$C$ a generic positive constant.
   We will also use $C_\epsilon$, $K_\epsilon$, $R_\epsilon$, etc.\ to denote constants that can depend on~$\epsilon$.

   The following proposition may be viewed as a first conditional stability result: 

\begin{prop}[Local stability for the lateral inverse problem]\label{prop_unilateral}
   Let us assume that
   \[
\Vert u^i\Vert_{H^2(Q_{\ell}(p_i))}\leq M,
\ \ \Vert p_i\Vert_{H^2(0,T)}\leq M
\ \text{ and } \ |p_i(t)|\leq 1-\delta \ \ \forall t\in[0,T]
   \]
for some~$\delta \in (0,1)$.
   Also, let us assume that~$0 < \epsilon<\bar{t}< T$ and
   \[
0\le D:= \Vert \alpha^1-\alpha^2\Vert_{L^2(0,T)}
+ \Vert \beta^1-\beta^2\Vert_{L^2(0,T)}<k<1 .
   \]
   Then there exist~$R_{\epsilon}$, $R_0>0$ and~$\mu_{\epsilon}\in(0,1)$ such that
   \begin{equation}\label{e_13a}
\Vert p_1-p_2\Vert_{L^{\infty}(\epsilon, T)}\\
\le \frac{R_{\epsilon}}{ \left( \log \frac{1}{k}\right)^{\mu_{\epsilon}}} 
+ R_0\vert p_1(\bar{t}) - p_2(\bar{t})\vert.
   \end{equation}
\end{prop}}

\

\noindent
\textbf{Proof:}  
   
   For instance, let us assume that~$p_1(t)\le p_2(t)$ for~$t\in(t_0,t_1)\subset[\epsilon,T]$ and  set~$h:=p_1-p_2$.
   Then, for all~$t\in(t_0,t_1)$ one has
   \begin{align*}
h'(t)&=u^1(p_1(t),t)-u^2(p_2(t),t)
\\[2mm]
& =\left[u^1(p_1(t),t)-u^2(p_1(t),t)\right]+\left[u^2(p_1(t),t)-u^2(p_2(t),t)
\right]
\\[2mm]
&\leq\frac{K_{\epsilon}}{\left(\log{\frac{1}{D}}\right)^{\theta_{\epsilon}}}
+2Mh(t) .
   \end{align*}
   Here, we have applied Lemma~\ref{lemma_key} to~$u^1-u^2$ in combination with 
the {\it Mean Value Theorem} to~$u^2(\cdot, t)$. 

   {\black It is clear that we can get a similar estimate in the whole interval~$[\epsilon, T]$:
   \[
h'(t)\leq\frac{K_{\epsilon}}{\left(\log{\frac{1}{D}}\right)
^{\theta_{\epsilon}}}+C\,|h(t)|
   \]
and, consequently
   \[
\frac{1}{2}\frac{d}{dt}|h(t)|^2\leq C\,|h(t)|^2
+ \frac{K_{\epsilon}}{\left(\log{\frac{1}{D}}\right)^{2\theta_{\epsilon}}} .
   \]}

   Hence, from~Gronwall's Lemma, we finally see that  
   \[
|h(t)|^2\leq \frac{K_{\epsilon}}{\left(\log{\frac{1}{D}}\right)^{2\theta_{\epsilon}}}
+ C \, |h(\bar{t})|^2 \quad\forall\, t\in[\epsilon, T]
   \]
and this leads to~\eqref{e_13a} 
\fin

\
 
\begin{remark}\label{cor_stab}
   Let the assumptions in~Proposition~\ref{prop_unilateral} be satisfied.
   Also, suppose that
   \[
\Vert u^i\Vert_{W^{2,\infty}(Q_{\ell}(p_i))}\leq M .
   \]
   Then, it can be ensured that for every~$\epsilon>0$ there exist~$K_{\epsilon}$, $K_0$ and~$\theta_{\epsilon}\in (0,1)$ such that
   \begin{equation}\label{e_13b}
\Vert p'_1-p'_2\Vert_{L^{\infty}(\epsilon, T)}
+ \Vert p''_1-p''_2\Vert_{L^{\infty}(\epsilon, T)}
\leq \dfrac{K_{\epsilon}}{\left(\log{\frac{1}{D}}\right)^{\theta_{\epsilon}}}
+K_0 \, |p_1(\bar{t})-p_2(\bar{t})|.
   \end{equation}
   Indeed, when~$p_1(t) \leq p_2(t)$, we have
   \[
p_1'(t)\!-\!p_2'(t)\!=\!u^1(p_1(t),t)\!-\!u^2(p_2(t),t)
\!=\! \left(u^1(p_1(t),t)\!-\!u^2(p_1(t),t)\right)\!+\!\left(u^2(p_1(t),t)\!-\!u^2(p_2(t),t)
\right)
   \]
and a similar identity holds when~$p_1(t)>p_2(t)$.
   Consequently, using Lemma~\ref{lemma_key} and the Mean Value Theorem, we can estimate~$\Vert p'_1-p'_2\Vert_{L^{\infty}(0, T)}$.
   On the other hand, 
   \begin{align*}
p_1''(t)-p_2''(t)=&\left(u^1_x(p_1(t),t)\,p_1'(t)+u^1_t(p_1(t),t)\right)
- \left(u_x^2(p_2(t),t)\,p_2'(t)+u^2_t(p_2(t),t)\right)\\[2mm]
=& \ u^1_x(p_1(t),t)(p_1'(t)-p_2'(t))+\left(u^1_x(p_1(t),t)-u^2_x(p_1(t),t)
\right)p_2'(t)
\\[2mm]
&+ \left(u_x^2(p_1(t),t)-u^2_x(p_2(t),t)\right)p_2'(t)
+ \left(u^1_t(p_1(t),t)-u^2_t(p_1(t),t)\right)
\\[2mm]
&+\left(u^2_t(p_1(t),t)-u^2_t(p_2(t),t)\right)
   \end{align*}
and, once more, a similar identity holds when~$p_1(t)>p_2(t)$. 
   Using Lemma~\ref{lemma_key}, the interpolation inequality, the Mean Value Theorem and the bounds
   \[
\Vert u^i_{xx}\Vert_{L^{\infty}}\leq M\quad\text{and}\quad\Vert u^i_{xt}\Vert_{L^{\infty}}
\leq M\quad\text{for}\quad i=1,2,
   \]
we can estimate~$\Vert p''_1-p''_2\Vert_{L^{\infty}(\epsilon, T)}$ and find~\eqref{e_13b}.
\fin
\end{remark}

\


\begin{cor}\label{cor_3.3}
   Under the assumptions in~Proposition~\ref{prop_unilateral}, if~$0<\overline{t}<T$ 
and~$\alpha^1\equiv\alpha^2$ and~$\beta^1\equiv\beta^2$ in~$(0,T)$, there exists a constant~$R_0>0$ such that 
   \begin{equation}\label{eq.cor33}
\Vert p_1-p_2\Vert_{L^{\infty}(0,T)}
\leq R_0 \, |p_1(\overline{t})-p_2(\overline{t})|,
   \end{equation}
where~$R_0$ is independent of~$\overline{t}$.
\end{cor}

\

\noindent
\textbf{Proof:}

   We can argue as in the proof of~Proposition~\ref{prop_unilateral}. 
   Thus, for every~$\epsilon>0$ and every small~$k>0$, we obtain 
   \[
\Vert p_1-p_2\Vert_{L^{\infty}(\epsilon, T)}\\
\le \frac{R_{\epsilon}}{ \left( \log \frac{1}{k}\right)^{\mu_{\epsilon}}} 
+ R_0\vert p_1(\overline{t}) - p_2(\overline{t})\vert.
   \]
   Then, taking~$k\to 0$, we see that 
   \[
\Vert p_1-p_2\Vert_{L^{\infty}(\epsilon, T)}\\
\le R_0\vert p_1(\overline{t}) - p_2(\overline{t})\vert.
   \]
   Finally, taking~$\epsilon \to 0$, we arrive at~\eqref{eq.cor33}. 
\fin

\


\begin{cor}[Lateral uniqueness]\label{cor_3.4}
   In addition to the assumptions in~Corollary~\ref{cor_3.3}, let us assume that~$p_1(\overline{t}) = p_2(\overline{t})$ for  some~$\overline{t}\in(0,T)$. 
   Then,
   \[
p_1\equiv p_2\quad\text{in}\quad(0,T)\quad\text{and}\quad u^1\equiv u^2
\quad\text{in}\quad Q_{\ell}(p).
   \]
\end{cor}

\


\begin{cor}\label{cor_3.5}
   In addition to the assumptions in~Proposition~\ref{prop_unilateral}, let us assume that~$p_1(\overline{t}) = p_2(\overline{t})$ for some $\overline{t}$ with~$0<\overline{t}<T$.
   Then, for any small~$\epsilon > 0$ there exist~$K_{\epsilon}>0$ and~$\theta_{\epsilon}\in (0,1)$ such that 
   \begin{equation}\label{e_15}
\Vert p_1-p_2\Vert_{L^{\infty}(\epsilon,T)}
\leq \frac{K_{\epsilon}}{\left(\log{\frac{1}{D}}\right)^{\theta_{\epsilon}}}.
   \end{equation}
\end{cor}

\
 
\begin{remark}
   If we do not assume that~$p_1(\overline{t})=p_2(\overline{t})$ for some~$\overline{t}\in(0,T)$, then we do not have uniqueness in general. 
   Indeed, the following particular functions furnish a counter-example:
   
\begin{itemize}

\item Let us set 
   \begin{align*}
&\begin{cases}
\ \va^1(x,t)=e^{-a^2t}\sin{ax}+A\quad\text{in}\quad Q_\ell:= (0,\ell)\times(0,T),
                              \\
\ \va^2(x,t)=e^{-a^2t}\sin{ax}+A\quad\text{in}\quad Q_L:=(0,L)\times(0,T),
\end{cases}
   \end{align*}
where~$a=n\pi/\ell=k\pi/L$ with~$0<\ell<L$ and~$n\neq k$ and~$A>0$.
   Then~$\va^i=\va^i(x,t)$ satisfies the heat equation for~$i=1,2$ and we observe that~$\va^1(\ell,t)\equiv\va^2(L,t)=A$.
   
\item Let us take
   \[
p_1(t)=\ell\quad\text{and}\quad p_2(t)=L\quad\text{in}\quad(0,T)
   \]
and set~$u^i:=-\va^i_x/\va^i$ for~$i=1,2$.
   Then the functions  $u^1$ and~$u^2$ satisfy the Burgers' equation respectively 
in~$Q_\ell$ and~$Q_L$ and do not coincide (see for example~\cite[Section 2.1.2]{ACDFY}).
   However the associated $\alpha^i$ and $\beta^i$ coincide.
\fin

\end{itemize}

\end{remark}

\

\section{Global estimates and uniqueness}\label{sec_global}

   In this section we will present global stability and uniqueness results for the inverse problem formulated in Section~\ref{sec_introduction} in the whole domain~$Q(p)$.

\

\begin{teo}[Conditional stability]\label{con_stability}
    Let~$(w_1,p^1)$ and~$(w_2,p^2)$ be the solutions to~\eqref{e_1} respectively corresponding 
to the data~$u_0$, $\alpha$, $m$, $q_0$, $q_1$, $\beta^i$ and~$\eta^i$ and set~$\beta^i(t)=w^i_x(-1,t)$ and~$\eta^i(t)= w_i(1, t)$ for 
$i=1,2$ and all~$t\in (0,T)$.
   {\black Assume that there exist constants~$\delta, \kappa \in (0,1)$ and~$M>0$ such that~$|p_i(t)|\leq 1-\delta$ for all~$t\in(0,T)$,
   \[
\Vert u^i\Vert_{W^{2,\infty}(Q_{\ell}(p_i))} \leq M , \ \ 
\Vert u^i\Vert_{W^{2,\infty}(Q_{r}(p_i))} \leq M \ \ (i=1,2) \ \text{ and }
\ 0 \leq \Vert \beta^1-\beta^2\Vert_{L^2(0,T)} < \kappa < 1 .
   \]
   Then, for every~$\epsilon>0$, there exist constants~$C_0, C_{\epsilon}>0$ and~$\theta_{\epsilon}\in (0,1)$ such that 
   \begin{equation}\label{e_17}
\Vert \eta^1-\eta^2\Vert_{L^{\infty}(\epsilon, T)}
\leq \frac{C_{\epsilon}}{\left|1+\log{\left|\log{\kappa}\right|}\right|
^{\theta_{\epsilon}}} + C_0|\eta^1(\overline{t})-\eta^2(\overline{t})|, 
   \end{equation}
for all~$\overline{t} \in [\epsilon, T)$.}
\end{teo}

\

\noindent
\textbf{Proof:}

   We will find estimates of~$ \Vert \eta^1-\eta^2\Vert_{L^{\infty}(\epsilon,T)}$ 
in terms of~$p_1-p_2$, $p'_1-p'_2$ and~$p''_1-p''_2$ and then we will use~\eqref{e_13a} and~\eqref{e_13b}.

   Let us set~$v^i=w_i|_{Q_r(p_i)}$, for~$i=1,2$.
   We introduce a change of  variables to move from~$Q_r(p_i)$ to~$(0,1)\times(0,T)$.
   For example, we set
   \[
\overline{v}^1(\overline{y},t) = v^1(x,t) \quad\text{with}\quad \overline{y} 
= \frac{x-p_1(t)}{1-p_1(t)}.
   \]
   Accordingly, we get the system
\begin{equation*}
\begin{cases}
\ \overline{v}^1_t -\dfrac{1}{(1-p_1(t))^2}\overline{v}^1
_{\overline{y}\overline{y}} 
- \dfrac{p_1'(t)(1-\overline{y})}{1-p_1(t)}\overline{v}^1_{\overline{y}} 
+ \dfrac{1}{1-p_1(t)}\overline{v}^1\overline{v}^1_{\overline{y}} = 0,    
&(\overline{y},t)\in (0,1)\times(0,T), 
\\[4mm]
\ \overline{v}^1(0,t) = p'_1(t), & t\in (0,T),
\\[2mm]
\ \overline{v}^1_{\overline{y}}(0,t)= (1-p_1(t))\,v^1_x(p_1(t),t), & t\in (0,T),
\\[2mm]
\ \overline{v}^1(1,t)=\eta^1(t), & t\in (0,T).
    \end{cases}
   \end{equation*}
   Furthermore, we note from~\eqref{eq_jump} that
   \[
(1-p_1(t))\,v^1_x(p_1(t),t)
= (1-p_1(t))\,u^1_x(p_1(t),t) + (1-p_1(t))\,p''_1(t)
   \]
for any~$t\in (0,T)$.

   A similar change of variables and a similar system hold for~$\overline{v}^2$ and~$p_2$. 

   Then, after some computations, for~$v:=\overline{v}^1-\overline{v}^2$ and~$p:=p_1-p_2$, we find that
   \begin{equation*}
\begin{cases}
\ v_t -\dfrac{1}{(1-p_1(t))^2}v_{\overline{y}\,\overline{y}} 
- \dfrac{p_1'(t)(1-\overline{y})}{1-p_1(t)}v_{\overline{y}} 
+ \dfrac{1}{1-p_1(t)}v = f, &   (\overline{y},t)\in (0,1)\times(0,T), 
\\[4mm]
\ v(0,t) = \widetilde{\alpha}(t), \quad  v_{\overline{y}}(0,t) 
= \widetilde{\beta}(t), & t\in (0,T),
    \end{cases}
   \end{equation*}
where~$\widetilde{\alpha}(t) = p_1'(t) -p_2'(t)$ and
   \begin{align*}
\widetilde{\beta}(t) &= (1-p_1(t))\,u^1_x(p_1(t),t) + (1-p_1(t))\,p''_1(t) \\
&\ \ -(1-p_2(t))\,u^2_x(p_2(t),t) - (1-p_2(t))\,p''_2(t)\\
& = \ (1-p_1(t))\,v^1_x(p_1(t),t) - (1-p_2(t))\,v^2_x(p_2(t),t). 
   \end{align*}
   
   One has
   \[
\Vert f\Vert_{L^{\infty}((0,1)\times(\epsilon, T))}
\leq C_1\left(\Vert p_1-p_2\Vert_{L^{\infty}(\epsilon, T)}+\Vert p'_1-p'_2\Vert 
_{L^{\infty}(\epsilon, T)}\right)
   \]
and
   \[
\Vert \widetilde{\alpha}\Vert_{L^{\infty}(\epsilon, T)}
+ \Vert \widetilde{\beta}\Vert_{L^{\infty}(\epsilon, T)}
\leq C_1\left(\Vert p_1-p_2\Vert_{L^{\infty}(\epsilon, T)}+\Vert p'_1-p'_2\Vert 
_{L^{\infty}(\epsilon, T)}
+ \Vert p''_1-p''_2\Vert_{L^{\infty}(\epsilon, T)}\right).
   \]
   Using now~\eqref{e_4b}, we see that
   \begin{equation}\label{e_18}
\Vert \eta^1-\eta^2\Vert_{L^{\infty}(\epsilon, T)}
\leq \frac{K_\epsilon}
{\left( \log [C_2/\Vert p_1-p_2\Vert_{W^{2,\infty}(\epsilon/2,T)}]\right)^{\theta_{\epsilon}}}
+ K_0|\eta^1(\overline{t})-\eta^2(\overline{t})|.
   \end{equation}
   Also, using~\eqref{e_13a} and~\eqref{e_13b}, we find 
   \begin{align}\label{e_19}
\Vert p_1-p_2\Vert_{w_{2,\infty}(\epsilon/2, T)}
&\leq\frac{K_{\epsilon}}{\left(\log{\frac{1}{\Vert \beta_1-\beta_2\Vert_{L^2(0,T)}}}
\right)^{\theta_{\epsilon}}}
+ K_0|p_1(\widetilde{t})-p_2(\widetilde{t})|   \nonumber  \\[4mm]
&\leq\frac{K'_{\epsilon}}{\left(\log{\frac{1}{\Vert \beta_1-\beta_2\Vert_{L^2(0,T)}}}
\right)^{\theta_{\epsilon}}}
   \end{align}
for~$\widetilde{t}$ close enough to~$0$ and~$\epsilon>0$ small enough.
   
   Hence, we obtain~\eqref{e_17} for some~$C_{\epsilon}$ and~$C_0$.

   This ends the proof.
\fin

\

{\black\begin{cor}[Global uniqueness]\label{cor_4.3}
   Let the assumptions in~Theorem~\ref{con_stability} be satisfied and let us assume that~$\beta^1 = \beta^2$ in~$(0,T)$.
   Then
   \begin{equation}
\eta^1=\eta^2 \quad\text{in}\quad (0,T). 
   \end{equation}
\end{cor}}

\

\noindent
\textbf{Proof:}

   Given an arbitrary~$\epsilon>0$, we take~$\bar{t}_\epsilon = 2\epsilon$. 
   For every~$\kappa>0$, using Theorem~\ref{con_stability} we can write that 
   \begin{equation*}
\Vert \eta^1-\eta^2\Vert_{L^{\infty}(\epsilon, T)}
\leq
\dfrac{C_{\epsilon}}{\left|1+\log{\left|\log{\kappa}\right|}
\right|^{\theta_{\epsilon}}}
+ C_0|\eta^1(\overline{t}_\epsilon)-\eta^2(\overline{t}_\epsilon)|.
   \end{equation*}

   Taking now~$\kappa\to 0$, we see that
   \begin{equation*}
\Vert \eta^1-\eta^2\Vert_{L^{\infty}(\epsilon, T)}
\leq C_0|\eta^1(\overline{t}_\epsilon)-\eta^2(\overline{t}_\epsilon)|.
   \end{equation*} 
   Finally, when~$\epsilon \to 0$, we see that
   \begin{equation*}
\Vert \eta^1-\eta^2\Vert_{L^{\infty}(0, T)}
\leq C_0|\eta^1(0)-\eta^2(0)|
   \end{equation*} 
and, since~$\eta^1(0) = \eta^2(0) = w_0(1)$, we deduce that  $\eta^1=\eta^2$ in~$(0,T)$. 
\fin


%

\section*{Acknowledgements}

   {\black 
   
 The first author was supported by the Grant PID2021-126813NB-I00 funded
by MCIN/AEI/10.1303\\ 9/501100011033 and by “ERDF A way of making Europe” and
by the grant IT1615-22 funded the Basque Government.
   
   The second and third authors were partially supported by MICINN, under Grant~PID2020-114976GB-I00. 
   
   The fourth author was supported by Grant-in-Aid for Scientific Research (A) 20H00117  and Grant-in-Aid for Challenging Research (Pioneering) 21K18142 of the Japan Society for the Promotion of Science.}


\begin{thebibliography}{999}
\addcontentsline{toc}{section}{References}


\bibitem{ACDFY} J. Apraiz, A. Doubova, E. Fern\'andez-Cara, M. Yamamoto, \emph{Some inverse problems for the Burgers equation and related systems}, Communications in Nonlinear Science and Numerical Simulation~\textbf{107} (2022), 106113.s


\bibitem{BK-1}
A.L.~Bukhgeim and M.V.~Klibanov,  \emph{Global uniqueness of a class of 
multidimensional inverse problems}, Soviet.\ Math. Doklady~\textbf{24}, (1981), 244--247.





\bibitem{DF2}
A.~Doubova and E.~Fern{\'a}ndez-Cara,  \emph{Some control results for simplified one-dimensional models of fluid-solid interaction}, Mathematical Models \& Methods in Applied Sciences~\textbf{15}, no.~5 (2005), 783--824. 














\bibitem{K-1}
M.V.~Klibanov, \emph{Inverse problems and Carleman estimates}, Inverse Problems~\textbf{8} (1992) 575--596.

\bibitem{K-2}
M.V.~Klibanov and A.A. Timonov, \emph{Carleman estimates for coefficient inverse 
problems and numerical applications}, Inverse and Ill-posed Problems Series, VSP, Utrecht, 2004.

\bibitem{LTT}
Y. Liu, T. Takahashi, and M. Tucsnak, \emph{Single input controllability of a simplified fluid-structure interaction model},  ESAIM: Control, Optimisation and Calculus of Variations~\textbf{19} (2013) 20--42.

\bibitem{MZ}
A. Munnier and E. Zuazua, Large time behavior for a simplified n-dimensional model of fluid–solid interaction, Communications in Partial Differential Equations~\textbf{30}, no.~3, (2005) 377--417.


\bibitem{VZ1} J.L.~V\'azquez and E. Zuazua, Lack of collision in a simplified 1-d model for fuid-solid inter\-action, Mathematical Models \& Methods in Applied Sciences~\textbf{16}, no.~5, (2006), 637--678.

\bibitem{VZ2} J.L.~V\'azquez and E. Zuazua, Large time behavior for a simplified 1D model of fluid-solid inter\-action,  Comm. Partial Differential Equations~\textbf{28} (2003), no.~9--10, 1705--1738. 
 
\bibitem{Ya} M. Yamamoto, \text{Carleman estimates for parabolic equations and  applications}, Inverse Problems~\textbf{25} (2009) 123013.


\end{thebibliography}
\end{document}